\documentclass[11pt]{article}
\usepackage{latexsym,amsmath,amssymb,epsfig,graphicx,rotating}
\usepackage[latin1]{inputenc}
\newtheorem{thm}{Theorem}[section]

\newtheorem{prop}[thm]{Proposition}
\newtheorem{defi}{Definition}
\newcommand{\Real}{\bf R}
\newcommand{\Nat}{\bf N}

\newenvironment{proof}{{\it Proof.}}
{\quad \hfill $\Box$}
\begin{document}

\title{\Large \bf A study on new computational local orders of convergence}
\author{{Miquel Grau-S\'{a}nchez$^{\rm a}$, Miquel Noguera$^{\rm a}$, \`{A}ngela Grau$^{\rm a}$, Jos\'{e} R. Herrero$^{\rm b}$} \\
                  \hfill    \\ [-0.6em]
 {\footnotesize $\,^{\rm a}${\it Technical University of Catalonia, Department of Applied Mathematics II}}\\ [-0.2em]
      {\footnotesize {\it Jordi Girona 1-3, Omega, 08034 Barcelona, Spain}} \\
     \hfill    \\ [-0.8em]
{\footnotesize $\,^{\rm b}${\it Technical University of Catalonia, Department of Computer Architecture}} \\ [-0.2em]
      {\footnotesize {\it Jordi Girona 1-3, C6, 08034 Barcelona, Spain}} \\ [0.2em]
       {\footnotesize  E-mail address: miquel.grau@upc.edu, miquel.noguera@upc.edu, angela.grau@upc.edu, josepr@ac.upc.edu}\\[0.2em]
       {\footnotesize  Research supported partially by the projects MTM2011-28636-C02-01}\\ [-0.2em]
       {\footnotesize and TIN2007-60625 of the Spanish Ministry of Science and Innovation}
      }

\date{}

\maketitle

%
%
%%%%%%%%%%%%%%%%%%%%%%%%%%%%%%%%%%%%%%%%%%%%%%%%%%%%%%%%%%%%%%%%%%%%%%%%%%%%%%%%%%%%%%%%%%%%%%%%%%%%%%
%%%%%%%%%%%%%%%%%%%%%%%%%%%%%%%%%%%%%%%%%%%%%%%%%%%%%%%%%%%%%%%%%%%%%%%%%%%%%%%%%%%%%%%%%%%%%%%%%%%%%%
%
%

%
% ----------------------------------------------------------------
\vspace{-9mm}
\begin{abstract}
Four new variants of the Computational Order of Convergence (COC) of a one-point iterative method with memory for solving nonlinear equations are presented. Furthermore, the way to approximate the new variants to the local order of convergence is analyzed. Three of the new definitions given here do not involve the unknown root. Numerical experiments using adaptive arithmetic with multiple precision and a stopping criteria are implemented without using any known root.
\end{abstract}

%
%
%%%%%%%%%%%%%%%%%%%%%%%%%%%%%%%%%%%%%%%%%%%%%%%%%%%%%%%%%%%%%%%%%%%%%%%%%%%%%%%%%%%%%%%%%%%%%%%%%%%%%%%%
%
%

\vspace{1.5mm}
{\small

\vspace{0.5mm}
}

%
%
%%%%%%%%%%%%%%%%%%%%%%%%%%%%%%%%%%%%%%%%%%%%%%%%%%%%%%%%%%%%%%%%%%%%%%%%%%%%%%%%%%%%%%%%%%%%%%%%%%%%%%%
%%%%%%%%%%%%%%%%%%%%%%%%%%%%%%%%%%%%%%%%%%%%%%%%%%%%%%%%%%%%%%%%%%%%%%%%%%%%%%%%%%%%%%%%%%%%%%%%%%%%%%%
%
%

\vspace{-2mm}
\section{\large Introduction}

\vspace{-2mm}
One-point iterative methods with memory for solving a nonlinear equation
$f(x) = 0 $, where $ f: I \subseteq \mathbb{R} \rightarrow \mathbb{R}$, and $I$ is a neighborhood of the root $\alpha$, usually consider a sequence $\{x_n\}_{n\in\Nat}$, defined by

\vspace{-6mm}
\begin{equation}\label{eqit}
  x_{n+1} = \phi \,(x_n; x_{n-1},\ldots, x_{n-j}) \, , \; \;\: n\ge 0,
\end{equation}

\noindent where $\phi$ is the iteration function. A sequence $\{x_n\}$ is said to converge to  $\alpha$ with \emph{local order of convergence}  $\rho  \in \Real$, $ \rho \ge 1$, if there exists the following limit

\vspace{-3mm}
\begin{equation} \label{order}
   \rho \,=\, \lim_{n \rightarrow\infty} \,  \frac{\log|e_{\scriptscriptstyle n+1}|}{\log |e_{\scriptscriptstyle n}|},
\end{equation}
where $\, e_k = x_k - \alpha$ is the error in the $k$th iterate (see \cite{Wal,Tor}). This limit $\rho$ is also equal to $R$-order defined in \cite{OrRh}. For one-point method with memory (\ref{eqit}) the error equation is:

\vspace{-2mm}
\begin{equation}\label{eqdif}
    e_{n+1} = C\,e_{n}^{\,\rho} \, \big( 1 \,+\, O (\,e_{n}^{\,\sigma}) \big),
\end{equation}
where $C$ is a real number, $0<\sigma<1$, and we will consider $\rho \ge (1+\sqrt{5})/2$. The nonzero constant $C$ is called the asymptotic error. The local order of convergence of an iterative method in a neighborhood of a root is the order of the corresponding sequence. If it is $\rho$, then the method approximately multiplies by $\rho$ the number of correct decimals after each iteration. That is, from (\ref{order}) we get $\log_{10} |e_{\scriptscriptstyle n+1}| \,\approx \, \rho \: \log_{10} |e_{\scriptscriptstyle n}|$, for $n$ large enough.

\vspace{1mm}
In the next sections the way to approximate four new variants of the local order of convergence is analyzed and numerical experiments using adaptive arithmetic with multiple precision and a stopping criteria are implemented without using any known root for three of the four techniques.

%%%%%%%%%%%%%%%%%%%%%%%%%%%%%%%%%%%%%%%%%%%%%%%%%%%%%%%%%%%%%%%%%%%%%%%%%%%%%%%%%%%%%%%%%%%%%%%%%%%%%
\section{\large Definitions and first result}

Next, we give the definitions of Computational Local Order of Convergence (CLOC) that is a variant of COC (\cite{WeFe}, 2000), Approximated Local Computational Order of Convergence (ACLOC), Extrapolated Local Computational Order of Convergence (ECLOC) and Petkovi\'{c} Local Computational Order of Convergence (PCLOC). These three last concepts are variants of ACOC, ECOC \cite{gng} and PCOC \cite{Petk} respectively. After the work of Weerakoon and Fernando \cite{WeFe}, many other authors have considered the COC in their research (see \cite{gn0}-\cite{gng} and references therein). In all those papers the COC is used to test numerically the order of convergence of the methods presented. Considering (\ref{order}) we provide a new parameter with lower cost than COC:

\vspace{-1mm}
\begin{defi}\hspace{-3.5mm}.  The computational local order of convergence (CLOC) of a sequence $\{x_n\}_{n\ge 0}$
is defined by

\vspace{-4mm}
\begin{equation}\label{CLOC}
     {\overline{\lambda}_n} = {\displaystyle \frac{\log |e_n|}{\log \left| e_{n-1} \right|},}
\end{equation}
where $x_{n-1}$ and $x_{n}$ are two consecutive iterations near the root $\alpha$ and $e_n =  x_n - \alpha$.
\end{defi}

\vspace{-1mm}
Notice that the last definition has lower cost because we use the logarithm function applied to only one variable, say $e_n$, instead of a quotient such as  $e_n/ e_{n-1}$ which is used in \cite{WeFe}.

\vspace{1mm}
The main drawback of COC and CLOC is that they involve the exact root $\alpha$, which in a real situation it is not known a priori. To avoid this, we introduce three variants of CLOC that do not use the exact root. Firstly, we give a new parameter considering three consecutive points:

\vspace{-2mm}
\begin{defi}\hspace{-2.5mm}. The approximated computational local order of convergence (ACLOC) of a sequence $\{x_n\}_{n\ge 0}$
is defined by

\vspace{-4mm}
\begin{equation}\label{ACLOC} \widehat{\lambda}_n =
 \frac{\log \,\left|\, \widehat{ e}_{n} \right|}{\log \, \left|\, \widehat{ e}_{n-1} \right|} \, ,
\end{equation}

\vspace{-3mm}
\noindent where $\, \widehat{ e}_n =  \, x_n - x_{n-1} \,$.
\end{defi}

\vspace{1mm}
Secondly, in order to avoid the requirement of the knowledge of the exact root $\alpha$, we consider three consecutive iterates $ x_n,\,x_{n-1},\,x_{n-2}$, and using Aitken's extrapolation we give the following approximation of $\alpha$

\vspace{-6mm}
\begin{equation}\label{xtsn}
{\widetilde \alpha}_n \,=\, x_n - \, \frac{ \left(\delta\, x_{n-1}\right)^2}{\delta^{\,2} \, x_{n-2}} \, , \quad n \ge 2,
\end{equation}

where $\delta$ is the forward difference operator, $ \delta x_k = x_{k+1} - x_k$ and (\ref{xtsn}) is the $\delta^2$-Aitken procedure \cite{Aitk}. Then, we can define a new approximation for the error $
{\widetilde e_n} =  x_n -{\widetilde \alpha}_n $ and a new computational order of convergence:

\vspace{-2mm}
\begin{defi}\hspace{-2.5mm}. The extrapolated computational local order of convergence (ECLOC) of a sequence $\{x_n\}_{n\ge 0}$ is defined by

\vspace{-4mm}
\begin{equation}\label{ECLOC} {\widetilde \lambda_n} =
 \frac{\log \,\left| \, {\widetilde e}_{n} \right|}{\log \, \left| \, {\widetilde e}_{n-1} \right| } \, ,
\end{equation}

\vspace{-1mm}
where ${\tilde e}_n =  x_n - {\tilde \alpha}_n  \,$ and $\,{\tilde \alpha}_n $ is given by (\ref{xtsn}).
\end{defi}

\vspace{1mm}
\noindent Finally, another way to avoid formulae involving the exact root $\alpha$ consists in using the values of two consecutive iterates. That is, from $ f(x_n)$ and $f(x_{n-1})$  the new computational order of convergence is:

\vspace{-2mm}
\begin{defi}\hspace{-2.5mm}. The Petkovi\'{c} computational local order of convergence (PCLOC) of a sequence $\{x_n\}_{n\ge 0}$ is defined by

\vspace{-4mm}
\begin{equation}\label{PCLOC} {\breve{ \lambda}_n} =
 \frac{\log \,\left| \, f(x_{n}) \right|}{\log \, \left| \, f(x_{n-1}) \right| } \, .
\end{equation}
\end{defi}

\vspace{-2mm}
This last parameter PCLOC is defined in honor of Petkovi\'{c} who in \cite{Petk,DPP} consider, in analogy of COC, the following value

\vspace{-6mm}
$$
  \breve{\rho}_n =  \dfrac {\log |f({x}_{\scriptstyle{n+1}}) / f(x_{\scriptstyle{n}})| }{\log |f(x_{\scriptstyle{n}}) / f(x_{\scriptstyle{n-1}})| } .
$$

As we show below, for all sequence $\{x_n\}$ converging to $\alpha$, with starting points $x_{-j},\ldots,x_{-1},x_0$ close enough to $\alpha$,  the values of ${\overline \lambda_n}$, ${\widehat \lambda_n}$, ${\widetilde \lambda_n}$ and $\breve{ \lambda}_n$   converge to  $\rho$, when $n \to \infty$.

\vspace{1mm}
There exist  numerical problems where a huge number of significant digits of the solution is needed. Such problems require the use of methods with a high order of convergence together with adequate arithmetics. We compute with a multiple precision arithmetic or symbolic manipulators, as Maple, that allow us to work with an adaptive arithmetic, that is, to update the length of the mantissa in each iteration by means of the formula

\vspace{-4mm}
\begin{equation}\label{eqdig1}
    {\tt Digits :=} \; \left[ \, \rho \times \left( - \log | \,e_n | +\, 2 \right) \right] \, ,
\end{equation}
where $\rho$ is the local order of convergence of the method and $[ x]$ denotes the integer part of  $ x $. Notice that the length of the mantissa is increased approximately by the order of convergence $\rho$.
In our experience, in order to guarantee all the significant  digits required we add $2$ units to $- \log | \,e_n |$.

%%%%%%%%%%%%%%%%%%%%%%%%%%%%%

\vspace{1mm}
Our first aim is to express $\,e_{n}$ as a function of $\,e_{n+1}$. In a first approximation, from (\ref{eqdif}), if we consider
$\, e_{n+1} \,=\, C\, e^{\,\rho}_{n}$ then we get $\,e_{n} \,=\, C^{\,- 1/\rho}\: e^{1/\rho}_{n+1}$. Substituting this result in the second term of the right side of (\ref{eqdif}) we obtain $\:
e_{n+1} \,=\, C\,e_{n}^{\,\rho} \, \big( 1 \,+\, O (\,e_{n+1}^{\,\sigma / \rho}) \big)\,
$, and

\vspace{-6mm}
\begin{eqnarray*}
 e_{n}^{\,\rho} &=& C^{\,-1}\, e_{n+1} \, \big( 1 \,+\, O (\,e_{n+1}^{\,\sigma / \rho}) \big).
\end{eqnarray*}

\vspace{-2mm}
Therefore, expressing $e_{n}$ in terms of $e_{n+1}$ we can state the following proposition:

%%%%%%%%%%%%%%%%%%%%%%%%%%%%%%%%%%%%%%%%%%%%%%%%%%%%%%%%%%%%%%%%%%%%%%%%%%%%%%%%%%%%%%%%%%%%%%%%%%%%%
%%%%%%%%%%%%%%%%%%%%%%%%%%%%%%%%%%%  Proposition 1.1   %%%%%%%%%%%%%%%%%%%%%%%%%%%%%%%%%%%%%%%%%%%%%%

\vspace{-1mm}
\begin{prop} \label{prop0}
Considering true the hypothesis in (\ref{eqdif}) we have

\vspace{-3mm}
\begin{equation}\label{enen1}
    e_n \:=\: C^{\,-1/\rho}\, e^{\, 1/\rho}_{n+1} \, \big( 1 \,+\, O (\,e_{n+1}^{\,\sigma / \rho}) \big).
\end{equation}
\end{prop}
%%%%%%%%%%%%%%%%%%%%%%%%%%%%%%%%%%%%%%%%%%%%%%%%%%%%%%%%%%%%%%%%%%%%%%%%%%%%%%%%%%%%%%%%%%%%%%%%%%%%%

%%%%%%%%%%%%%%%%%%%%%%%%%%%%%%%%%%%%%%%%%%%%%%%%%%%%%%%%%%%%%%%%%%%%%%%%%%%%%%%%%%%%%%%%%%%%%%%%%%%%%
%%%%%%%%%%%%%%%%%%%%%%%%%%%%%%%%%%%   CLOC   %%%%%%%%%%%%%%%%%%%%%%%%%%%%%%%%%%%%%%%%%%%%%%%%%%%%%%%%

\vspace{-5mm}
\section{\large Computational Local Order of Convergence (CLOC)}

\vspace{-2mm}
A relationship between  $\overline{\lambda}_n$ and $\rho$ is derived. In fact, we prove that $\, {\overline{\lambda}_n}$ converges to $\rho\,$ when $\, e_{n-1}\rightarrow 0$. That is $\,{\overline{\lambda}_n} \approx \rho$, in the sense that $\,{\displaystyle \lim_{n\rightarrow\infty } \frac{\overline{\lambda}_n}{\rho} = 1 }$.

%%%%%%%%%%%%%%%%%%%%%%%%%%%%%%%%%%%%%%%%%%%%%%%%%%%%%%%%%%%%%%%%%%%%%%%%%%%%%%%%%%%%%%%%%%%%%%%%%%%%%
%%%%%%%%%%%%%%%%%%%%%%%%%%%%%%%%%%  Proposition 2.1   %%%%%%%%%%%%%%%%%%%%%%%%%%%%%%%%%%%%%%%%%%%%%%%

\vspace{-1mm}
\begin{prop} If $\overline{\lambda}_n$ is the  CLOC defined in (\ref{CLOC}) and $\rho$ is the order of convergence, then

\vspace{-2mm}
\begin{equation} \label{prop1}
{\overline{\lambda}_n} \,=\, \rho \left( \, 1 + \: O \left( \frac{\log |C|}{\rho \, \log |e_{n-1}|}\right) \right),
\end{equation}

\vspace{-3mm}
where $C$ is given in (\ref{eqdif}).
\end{prop}
%%%%%%%%%%%%%%%%%%%%%%%%%%%%%%%%%%%%%%%%%%%%%%%%%%%%%%%%%%%%%%%%%%%%%%%%%%%%%%%%%%%%%%%%%%%%%%%%%%%%%
%%%%%%%%%%%%%%%%%%%%%%%%%%%%%%%%%%  Proof Proposition 2.1   %%%%%%%%%%%%%%%%%%%%%%%%%%%%%%%%%%%%%%%%%

\vspace{-1mm}
\begin{proof}
To prove (\ref{prop1}) we express ${\overline{\lambda}_n}$ in terms of $e_{n-1}$. Taking into account (\ref{eqdif}) we take

\vspace{-5mm}
\begin{eqnarray} \nonumber
 \log |e_{n}| &=&\log \left| C\,e^{\rho}_{n-1} \left( 1 +\,O( e^{\,\sigma}_{n-1})\right) \right|\\ [0.3em]\label{lnen}
              &=&\rho \log |e_{n-1}| + \log |C| + O(e^{\,\sigma}_{n-1}).
\end{eqnarray}

\vspace{-4mm}
Then

\vspace{-2mm}
\begin{equation} \nonumber
 {\overline{\lambda}_n} \,=\, \frac{\log |e_{n}|}{\log |e_{n-1}|} \,=\, \frac{\rho \log |e_{n-1}| + \log |C| + O(e^{\,\sigma}_{n-1})}{\log |e_{n-1}|},
\end{equation}

\vspace{1mm}
and we obtain the assertion of the proposition. %\hspace{3mm} $\Box$
\end{proof}
%%%%%%%%%%%%%%%%%%%%%%%%%%%%%%%%%%%%%%%%%%%%%%%%%%%%%%%%%%%%%%%%%%%%%%%%%%%%%%%%%%%%%%%%%%%%%%%%%%%%%%%%

\vspace{1mm}
Notice that for the calculus of the  CLOC (\ref{CLOC}) and for updating the adaptive arithmetic process (\ref{eqdig1}) it is necessary to know the exact root $\alpha $. In this case the following stopping criterion is applied:

\vspace{-3mm}
\begin{equation}\label{stop1}
     | e_n | = |x_n - \alpha|< 10^{-\eta},
\end{equation}

where $\eta$ is the maximum number of correct decimals and $ 10^{-\eta}$ is the required accuracy.

%%%%%%%%%%%%%%%%%%%%%%%%%%%%%%%%%%%%%%%%%%%%%%%%%%%%%%%%%%%%%%%%%%%%%%%%%%%%%%%%%%%%%%%%%%%%%%%%%%%%%%%%

%%%%%%%%%%%%%%%%%%%%%%%%%%%%%%%%%%%%%%%%%%%%%%%%%%%%%%%%%%%%%%%%%%%%%%%%%%%%%%%%%%%%%%%%%%%%%%%%%%%%%%%%
%%%%%%%%%%%%%%%%%%%%%%%%%%%%%%%%%%%   ACLOC   %%%%%%%%%%%%%%%%%%%%%%%%%%%%%%%%%%%%%%%%%%%%%%%%%%%%%%%%%%

\vspace{-3mm}
\section{\large Approximated Computational Local Order of Convergence (ACLOC)}

\vspace{-2mm}
A relationship between $\widehat{\lambda}_n$ and $\rho$ is obtained. A new technique to update the number of significant digits in an adaptive multi-precision arithmetic is given and a new stopping criterion is suggested.

%%%%%%%%%%%%%%%%%%%%%%%%%%%%%%%%%%%%%%%%%%%%%%%%%%%%%%%%%%%%%%%%%%%%%%%%%%%%%%%%%%%%%%%%%%%%%%%%%%%%%%%%
%%%%%%%%%%%%%%%%%%%%%%%%%%%%%%%%%%%% Proposition 4.1   %%%%%%%%%%%%%%%%%%%%%%%%%%%%%%%%%%%%%%%%%%%%%%%%%

\vspace{-1mm}
\begin{prop}  \label{P41} Let $\widehat{\lambda}_n$ be  the  ACLOC defined in (\ref{ACLOC}) and $\rho$ the order of convergence, then

\vspace{-2mm}
\begin{equation} \label{prop3}
\widehat{\lambda}_n \,=\, \rho \left( \, 1 + \: O \left( \frac{\log |C|}{\rho \, \log |e_{n-2}|}\right) \right),
\end{equation}

\vspace{-2mm}
where $C$ is given in (\ref{eqdif}).
\end{prop}
%%%%%%%%%%%%%%%%%%%%%%%%%%%%%%%%%%%%%%%%%%%%%%%%%%%%%%%%%%%%%%%%%%%%%%%%%%%%%%%%%%%%%%%%%%%%%%%%%%%%%%%%
%%%%%%%%%%%%%%%%%%%%%%%%%%%%%%%%%%%% Proof Proposition 4.1   %%%%%%%%%%%%%%%%%%%%%%%%%%%%%%%%%%%%%%%%%%%

\vspace{-1mm}
\begin{proof}
From the following expression:

\vspace{-3mm}
\begin{equation}\label{logeb}
     \log |\widehat{e}_n| \,=\, \log |{e}_n - {e}_{n-1}| \,=\, \log |{e}_{n-1}| + O(|{e}_{n}/{e}_{n-1}|),
\end{equation}

\vspace{-1mm}
and (\ref{ACLOC}), we get

\vspace{-4mm}
$$ \widehat{\lambda}_n \,=\,
 \frac{\log \,\left|\, \widehat{ e}_{n} \right|}{\log \, \left|\, \widehat{ e}_{n-1} \right|} \,=\,
  \frac{\log \,\left|\,{ e}_{n-1}\right| + O(\left| e_n/e_{n-1} \right|)}{\log \, \left|\,{ e}_{n-2}\right|
  + O(\left| e_{n-1}/e_{n-2} \right|)}.
$$

\vspace{-1mm}
Applying (\ref{lnen}) we obtain

\vspace{-5mm}
\begin{eqnarray*}
\widehat{\lambda}_n &=& \frac{\rho \log |e_{n-2}| + \log |C| + O(e^{\,\nu}_{n-2})}{\log \,
\left|\,{ e}_{n-2} \right| + O(\left| e_{n-1}/e_{n-2} \right|)},
\end{eqnarray*}

\vspace{-1mm}
 where $\, \nu = \min \{\rho(\rho-1), \sigma \}$ and the proof is complete.
\end{proof}
%%%%%%%%%%%%%%%%%%%%%%%%%%%%%%%%%%%%%%%%%%%%%%%%%%%%%%%%%%%%%%%%%%%%%%%%%%%%%%%%%%%%%%%%%%%%%%%%%%%%%%%%%

\vspace{1mm}
Observe that from (\ref{prop1}) and (\ref{prop3}) the expressions of $\overline{\lambda}_n$ and
$\widehat{\lambda}_{n+1}$ are identical. That is, if we approximate the theoretical value of the local order
$\rho$ by the computational values $\overline{\lambda}_n$ and $\widehat{\lambda}_{n}$:

\vspace{-2mm}
$$
\rho  = \overline{\lambda}_n \pm \Delta \overline{\lambda}_n , \quad
\rho  = \widehat{\lambda}_{n+1} \pm \Delta \widehat{\lambda}_{n+1} ,
$$

\vspace{-1mm}
then $ \Delta \widehat{\lambda}_{n+1} \approx \Delta \overline{\lambda}_n$. Our numerical experiments confirm this relation.
%%%%%%%%%%%%%%%%%%%%%%%%%%%%%%%%%%%%%%%%%%%%%%%%%%%%%%%%%%%%%%%%%%%%%%%%%%%%%%%%%%%%%%%%%%%%%%%%%%%%%%%%%

\vspace{2mm}
Repeating (\ref{eqdif}) twice we obtain

\vspace{-4mm}
$$
  e_n \,=\, C^{\,\rho + 1} \, e^{\rho^2}_{n-2} \big( 1 +\, O (e^{\sigma}_{n-2}) \big),
$$

\vspace{-1mm}
and now, we write $\widehat {e}_{n}/{\widehat e}_{n-1}$ in terms of $e_{n-2}$:

\vspace{-6mm}
\begin{eqnarray} \nonumber
  \frac{{\widehat e}_{n}}{{\widehat e}_{n-1}} \,=\, \frac{e_{n} - e_{n-1}}{e_{n-1} -e_{n-2}} & = &
  \frac{\,C^{\, \rho +1}   \, e^{\,\rho^2}_{n-2} + O (e^{\rho^2 + \sigma}_{n-2}) - C \, e^{\,\rho }_{n-2} + O ( e^{\,\rho + \sigma }_{n-2})}{C\,
  e^{\rho}_{n-2} + O ( e^{\,\rho + \sigma}_{n-2} ) -  e_{n-2} } \\ [0.6em] \label{errc11}
  & = & C\, e^{\,\rho -1}_{n-2} \left( 1 +  O \left( e^{\,\tau}_{n-2} \right) \right),
\end{eqnarray}

\vspace{-1mm}
where $ \tau = \min \,\{ \rho-1, \, \sigma\}$. Moreover, we get

\vspace{-2mm}
\begin{equation} \label{en2en}
e_{n-2} \:=\: C^{\,-(\rho+1)/ \rho^2} \: e^{\,1/\rho^2 }_{n} \, \big( 1 + O (e^{\,\sigma/\rho}_{n}) \big),
\end{equation}

\vspace{1mm}
since, from (\ref{enen1}), $e_{n-2} = C^{\,-1/ \rho} \: e^{\,1/\rho}_{n-1} \, \big( 1 + O (e^{\,\sigma/\rho}_{n-1}) \big)$ and $e_{n-1} = C^{\,-1/ \rho} \: e^{\,1/\rho}_{n} \, \big( 1 + O (e^{\,\sigma/\rho}_{n}) \big)$.

\vspace{2mm}
Substituting (\ref{en2en}) in (\ref{errc11}) we have the following proposition:

%%%%%%%%%%%%%%%%%%%%%%%%%%%%%%%%%%%%%%%%%%%%%%%%%%%%%%%%%%%%%%%%%%%%%%%%%%%%%%%%%%%%%%%%%%%%%%%%%%%%%
%%%%%%%%%%%%%%%%%%%%%%%%%%%%%%%%%%%%%%  Proposition 4.2   %%%%%%%%%%%%%%%%%%%%%%%%%%%%%%%%%%%%%%%%%%%

\vspace{-1mm}
\begin{prop} If we set  $e_n=x_n-\alpha$  and $ \widehat{ e}_n =  x_n -x_{n-1} $, then

\vspace{-2mm}
\begin{equation}\label{prop4}
e_n \approx \: C^{\,1 / (1-\rho)} \, \left( \frac{{\widehat e}_{n}}{{\widehat e}_{n-1}} \right)^{\rho^2 / (\rho -1)},
\end{equation}
where $\rho$ is the order of convergence and $C$ is given in (\ref{eqdif}).
\end{prop}
%%%%%%%%%%%%%%%%%%%%%%%%%%%%%%%%%%%%%%%%%%%%%%%%%%%%%%%%%%%%%%%%%%%%%%%%%%%%%%%%%%%%%%%%%%%%%%%%%%%%%

The result (\ref{prop4}) allows us to substitute the error term in (\ref{eqdig1}) by an expression that does not involve the exact root. Indeed, we implement the following adaptive multi-precision arithmetic scheme:

\vspace{-2mm}
\begin{equation}\label{eqdig3}
{\tt Digits } := \left[ \frac{\rho^3}{\rho-1} \times \left(- \log \left| \, \frac{{\widehat e}_{n}}{{\widehat e}_{n-1}} \right| \,+\, 2 \right) \,\right].
\end{equation}

Moreover, from (\ref{prop4}) we propose the following stopping criterion, instead of (\ref{stop1}):

\vspace{-2mm}
\begin{equation}\label{stop3}
     \left| \, \frac{{\widehat e}_{n}}{{\widehat e}_{n-1}} \right| <  10^{-\eta \: (\rho -1) / \rho^2} .
\end{equation}

%%%%%%%%%%%%%%%%%%%%%%%%%%%%%%%%%%%%%%%%%%%%%%%%%%%%%%%%%%%%%%%%%%%%%%%%%%%%%%%%%%%%%%%%%%%%%%%%%%%%%%
%%%%%%%%%%%%%%%%%%%%%%%%%%%%%%%%%%%   ECLOC   %%%%%%%%%%%%%%%%%%%%%%%%%%%%%%%%%%%%%%%%%%%%%%%%%%%%%%%%

\vspace{-3mm}
\section{\large Extrapolated Computational Local Order of Convergence (ECLOC)}

\vspace{-2mm}
We give a relationship between $\widetilde{\lambda}_n$ and  $\rho$, a new technique to update the number of significant digits in an adaptive multi-precision arithmetic and a new stopping criterion.

\vspace{1mm}
We start by deriving an expression of ECLOC as a function of the local order of convergence.
%%%%%%%%%%%%%%%%%%%%%%%%%%%%%%%%%%%%%%%%%%%%%%%%%%%%%%%%%%%%%%%%%%%%%%%%%%%%%%%%%%%%%%%%%%%%%%%%%%%%%%
%%%%%%%%%%%%%%%%%%%%%%%%%%%%%%%%%%%%%%  Proposition 5.1   %%%%%%%%%%%%%%%%%%%%%%%%%%%%%%%%%%%%%%%%%%%%

\vspace{-1mm}
\begin{prop} \label{P51} If $\,\widetilde{\lambda}_n$ is the ECLOC defined in (\ref{ECLOC}) and  $\rho \ge \phi$ is the order of convergence, then

\vspace{-3mm}
\begin{equation} \label{prop2}
\widetilde{\lambda}_n \approx \rho \left( \, 1 + \frac{\log |C| }{(2 \rho - 1)\, \log |e_{n-2}|}
    \right) ,
\end{equation}

\vspace{-1mm}
where $C$ is given in (\ref{eqdif}).
\end{prop}
%%%%%%%%%%%%%%%%%%%%%%%%%%%%%%%%%%%%%%%%%%%%%%%%%%%%%%%%%%%%%%%%%%%%%%%%%%%%%%%%%%%%%%%%%%%%%%%%%%%%%
%%%%%%%%%%%%%%%%%%%%%%%%%%%%%%%%%%%%%% Proof Proposition 5.1   %%%%%%%%%%%%%%%%%%%%%%%%%%%%%%%%%%%%%%

\vspace{-1mm}
\noindent \begin{proof}
Now, we write  $\widetilde{\lambda}_n$ in terms of $e_{n-2}$. To do that, we express $\,\log |\widetilde e_n|$ as a function of $e_{n-1}$ and $e_{n-2}$:

\vspace{-6mm}
\begin{eqnarray} \nonumber
\log |{\widetilde e}_n| &=& 2\,\log |\,e_n - e_{n-1}\,| - \log |\,e_n-2 e_{n-1} + e_{n-2}\,|\\ [0.3em] \nonumber
  &=& 2\, \log |\, e_{n-1}\,| + O( \left| \,e_n / e_{n-1}\,\right|) \, -\, \log |\,e_{n-2}\,| + O( \left| \,e_{n-1} / e_{n-2}\,\right|)\\ [0.3em]  \label{errt03}
 &=& \log \left|\, e^2_{n-1} / e_{n-2}\,\right| + O( \left| \,e_{n-1} / e_{n-2}\,\right|).
\end{eqnarray}

\vspace{-2mm}
We obtain

\vspace{-4mm}
$$
{\widetilde \lambda_n} =
 \frac{\log \,\left| \, {\widetilde e}_{n} \right|}{\log \, \left| \, {\widetilde e}_{n-1} \right| } \,=\,
\frac{\log \left|\, e^2_{n-1} / e_{n-2}\,\right| + O( \left| \,e_{n-1} / e_{n-2}\,\right|)}{\log \left|\, e^2_{n-2} / e_{n-3}\,\right| + O( \left| \,e_{n-2} / e_{n-3}\,\right|)}.
$$

From (\ref{eqdif}) we deduce $\, e^2_{n-1} = C^2 \, e^{2 \rho}_{n-2}\, \big( 1 + O(e^{\sigma}_{n-2}) \big)$, and taking into account (\ref{enen1}), we have $\, e_{n-3} = C^{\,-1/\rho} \, e^{1/\rho}_{n-2}\, \big( 1 + O(e^{\sigma / \rho}_{n-2}) \big)$. Next, we get

\vspace{-5mm}
\begin{eqnarray*} \nonumber
{\widetilde \lambda}_{n} &\approx & \frac{\log \left|\,C^2 \, e^{2 \rho - 1}_{n-2}\, \big( 1 + O(e^{\sigma}_{n-2}) \big)\, \right| }{\log \left|\,C^{\,1/\rho} \, e^{2- 1/\rho}_{n-2}\, \big( 1 + O(e^{\sigma / \rho}_{n-2}) \big)\,\right|}\\ [0.2em]
 &\approx & \rho \: \frac{( 2 \rho - 1) \log |e_{n-2}| + 2 \log |C| + O(e^{\sigma}_{n-2})}{( 2 \rho - 1) \log |e_{n-2}| + \log |C| + O(e^{\sigma / \rho}_{n-2})},
\end{eqnarray*}

\vspace{-3mm}
and the proof is complete.
\end{proof}
%%%%%%%%%%%%%%%%%%%%%%%%%%%%%%%%%%%%%%%%%%%%%%%%%%%%%%%%%%%%%%%%%%%%%%%%%%%%%%%%%%%%%%%%%%%%%%%%%%%%%
%%%%%%%%%%%%%%%%%%%%%%%%%%%%%%%%%%%%  Proposition 5.2   %%%%%%%%%%%%%%%%%%%%%%%%%%%%%%%%%%%%%%%%%%%%%

\begin{prop} \label{P52} Given $e_n=x_n-\alpha$ and  ${\widetilde e_n} =  x_n -{\widetilde \alpha_n} $, then

\vspace{-2mm}
\begin{equation}\label{bigeq}
e_n \approx C^{\:\beta} \: {\widetilde e}_n^{\; \, \rho^2 /\, (2\rho-1)} ,
\quad \mbox{where} \quad \beta = \frac{\rho-1}{2\rho-1}.
\end{equation}
\end{prop}
%%%%%%%%%%%%%%%%%%%%%%%%%%%%%%%%%%%%%%%%%%%%%%%%%%%%%%%%%%%%%%%%%%%%%%%%%%%%%%%%%%%%%%%%%%%%%%%%%%%%%
%%%%%%%%%%%%%%%%%%%%%%%%%%%%%%%%%%%%  Proof Proposition 5.2   %%%%%%%%%%%%%%%%%%%%%%%%%%%%%%%%%%%%%%%

\vspace{-1mm}
\noindent \begin{proof}
Taking into account $ \;
e_{n-2} = C^{\,- 1 / \rho}\,e^{\,1 / \rho}_{n-1} \left( 1 + O ( e_{n-1}^{\,\sigma / \rho}) \right),
$
we write ${\tilde e}_n$ in terms of $e_{n-1}$:

\vspace{-3mm}
\begin{eqnarray} \nonumber
  {\tilde e}_n \,=\, \frac{( e_{n} - e_{n-1})^2}{e_n- 2 e_{n-1} +e_{n-2}} & = & \frac{\,\big(C^{\, 2}
  \, e^{\,2\rho}_{n-1}  - 2 C \, e^{\,\rho + 1}_{n-1} \big) \big( 1 + O(e^{\sigma / \rho}_{n-1}) \big) + e^{\,2}_{n-1}}{
C\, e^{\rho}_{n-1} \big( 1 + O(e^{\sigma / \rho}_{n-1}) \big) - 2 e_{n-1} + \, C^{\,- 1/\rho} \: e^{1 / \rho}_{n-1} \left( 1 + O(e_{n-1}^{\,\sigma / \rho})
  \right)} \\ [0.3em]  \label{errt01}
  & = & C^{\, 1 / \rho}\, e^{(2\rho -1)/\rho}_{n-1} \left( 1 + O(e_{n-1}^{\,\tau / \rho}) \right),
\end{eqnarray}

\vspace{-2mm}
\noindent where $\, \tau = \min \{\rho-1, \sigma \}$. Now, from (\ref{errt01}) and $ \;
       e_{n-1} = C^{\,- 1 / \rho}\, e^{\,1 / \rho}_{n} \left( 1 + O ( e_n^{\, \sigma / \rho})\right),
$ we get

\vspace{-4mm}
\begin{eqnarray} \nonumber
  {\tilde e}_n &\hspace{-1mm}= &\hspace{-1mm} C^{\, 1/\rho} \left[ C^{\,- 1 / \rho}\,e^{\,1 / \rho}_{n} \left( 1 + O ( e_n^{\,\sigma / \rho}) \right)\right]^{(2\rho-1)/\rho} \cdot \left[ 1+ O \left( \left\{
    C^{\,- 1 / \rho}\, e^{\,1 / \rho}_{n} \left( 1 + O ( e_n^{\,\sigma / \rho}) \right) \right\}^{\tau /
    \rho} \right) \right] \\ [0.4em] \label{errt02}
     &\hspace{-1mm}=&\hspace{-1mm} C^{ \, (1-\rho ) / \rho^2} \: e^{\, ( 2 \rho -1 ) / \rho^2}_n \: \left( 1 + \, O (\,e^{\,  \tau / \rho^2}_n ) \right).
\end{eqnarray}

\vspace{-2mm}
\noindent From (\ref{errt02}),  we have
$ \;
 e^{\, ( 2 \rho -1 ) / \rho^2}_n \approx C^{ \, (\rho - 1 ) / \rho^2}  \:  {\tilde e}_n \,
$ from which the proof immediately follows.
\end{proof}
%%%%%%%%%%%%%%%%%%%%%%%%%%%%%%%%%%%%%%%%%%%%%%%%%%%%%%%%%%%%%%%%%%%%%%%%%%%%%%%%%%%%%%%%%%%%%%%%%%%%%

\vspace{2mm}
Notice that  (\ref{bigeq}) allows us to implement an iterative method (\ref{eqit}) with a multi-precision adaptive arithmetic. Instead of (\ref{eqdig1}) we now consider the expression:

\vspace{-2mm}
\begin{equation}\label{eqdig2}
{\tt Digits } := \left[ \frac{\rho^3}{2\rho-1} \times \left(- \log | \,\tilde e_n| \,+\, 2 \right) \,\right].
\end{equation}

\vspace{1mm}
In addition, as an alternative to  (\ref{stop1}), (\ref{bigeq}) provides the following stopping criterion

\vspace{-2mm}
\begin{equation}\label{stop2}
     | \, {\tilde e}_n  | <  10^{-\eta\: (2 \rho -1) / \rho^2} .
\end{equation}

%%%%%%%%%%%%%%%%%%%%%%%%%%%%%%%%%%%%%%%%%%%%%%%%%%%%%%%%%%%%%%%%%%%%%%%%%%%%%%%%%%%%%%%%%%%%%%%%%%%%%
%%%%%%%%%%%%%%%%%%%%%%%%%%%%%%%%%%%%%%   PCLOC   %%%%%%%%%%%%%%%%%%%%%%%%%%%%%%%%%%%%%%%%%%%%%%%%%%%%

\vspace{-3mm}
\section{\large Petkovi\'{c} Computational Local Order of Convergence (PCLOC)}

\vspace{-2mm}
In this section we provide a relationship between $\breve{\lambda}_n$ and $\rho$. In addition, we derive a new technique to update the number of significant digits in an adaptive multi-precision arithmetic and a new stopping criterion.

%%%%%%%%%%%%%%%%%%%%%%%%%%%%%%%%%%%%%%%%%%%%%%%%%%%%%%%%%%%%%%%%%%%%%%%%%%%%%%%%%%%%%%%%%%%%%%%%%%%%%
%%%%%%%%%%%%%%%%%%%%%%%%%%%%%%%%%%%%  Proposition 6.1   %%%%%%%%%%%%%%%%%%%%%%%%%%%%%%%%%%%%%%%%%%%%%

\vspace{-2mm}
\begin{prop} \label{P61} If $\breve{\lambda}_n$ is  the  PCLOC defined in (\ref{PCLOC}) and $\rho$ is the order of convergence, then

\vspace{-1mm}
\begin{equation}
\breve{\lambda}_n \,\approx\, \rho \left( \, 1 + \: O \left( \frac{\log | C\, \Gamma^{1-\rho} |}{\rho \, \log |e_{n-1}|}\right) \right),
\end{equation}

\vspace{-1mm}
where $C$ is given in (\ref{eqdif}) and $\Gamma = f'(\alpha)$.
\end{prop}
%%%%%%%%%%%%%%%%%%%%%%%%%%%%%%%%%%%%%%%%%%%%%%%%%%%%%%%%%%%%%%%%%%%%%%%%%%%%%%%%%%%%%%%%%%%%%%%%%%%%%
%%%%%%%%%%%%%%%%%%%%%%%%%%%%%%%%%%%%  Proof Proposition 6.1   %%%%%%%%%%%%%%%%%%%%%%%%%%%%%%%%%%%%%%%

\vspace{-1mm}
\begin{proof}
Setting
\begin{equation}\label{Tayf}
     f(x_k) = \Gamma \left( e_k + O (e^2_k) \right) = \Gamma e_k \left( 1 + O (e_k) \right), \; \; k = n-1, n ,
\end{equation}
and from (\ref{eqdif}),   $e_n = C e^{\rho}_{n-1} \left( 1 + O(e^{\sigma}_{n-1}\right)$,  we have

\vspace{-2mm}
$$ \breve{\lambda}_n \,=\,
 \frac{\log \left| f(x_{n}) \right|}{\log \left|\, f(x_{n-1}) \right|} \,=\,
  \frac{\rho \,\log \left| { e}_{n-1} \right| + \log |\Gamma C| +  O(\left| e^{\sigma}_{n-1} \right|)}{\log  \left| { e}_{n-1}\right| + \log |\Gamma| + O(\left| e_{n-1} \right|)}.
$$

\vspace{1mm}
Dividing the numerator and denominator in the right side of the preceding equation by $\log \, \left|\,{ e}_{n-1}\right| $ the proof is complete.
\end{proof}
%%%%%%%%%%%%%%%%%%%%%%%%%%%%%%%%%%%%%%%%%%%%%%%%%%%%%%%%%%%%%%%%%%%%%%%%%%%%%%%%%%%%%%%%%%%%%%%%%%%%%

%%%%%%%%%%%%%%%%%%%%%%%%%%%%%%%%%%%%%%%%%%%%%%%%%%%%%%%%%%%%%%%%%%%%%%%%%%%%%%%%%%%%%%%%%%%%%%%%%%%%%
%%%%%%%%%%%%%%%%%%%%%%%%%%%%%%%%%%%%%%  Proposition 6.2   %%%%%%%%%%%%%%%%%%%%%%%%%%%%%%%%%%%%%%%%%%%

\begin{prop} If we denote $Q_n=\dfrac{f(x_n)}{f(x_{n-1})}$,  then
\begin{equation}\label{prop5}
    e_n \approx C^{\,\frac{1}{1- \rho}} \; {Q}^{\,\frac{\rho}{\rho-1}}_n .
\end{equation}
\end{prop}
%%%%%%%%%%%%%%%%%%%%%%%%%%%%%%%%%%%%%%%%%%%%%%%%%%%%%%%%%%%%%%%%%%%%%%%%%%%%%%%%%%%%%%%%%%%%%%%%%%%%%
%%%%%%%%%%%%%%%%%%%%%%%%%%%%%%%%%%%%  Proof Proposition 6.2   %%%%%%%%%%%%%%%%%%%%%%%%%%%%%%%%%%%%%%%

\vspace{-2mm}
\noindent \begin{proof}
Taking into account (\ref{Tayf}) and (\ref{enen1}) we have  $\, Q_n \,=\, C^{\,\frac{1}{\rho}} \; {e}^{\,\frac{\rho-1}{\rho}}_n \left( 1 + O(e^{\sigma/\rho}_n) \right)$ and the proof immediately follows.
\end{proof}
%%%%%%%%%%%%%%%%%%%%%%%%%%%%%%%%%%%%%%%%%%%%%%%%%%%%%%%%%%%%%%%%%%%%%%%%%%%%%%%%%%%%%%%%%%%%%%%%%%%%%

\vspace{2mm}
The result (\ref{prop5}) allows us to substitute the error in (\ref{eqdig1}) by an expression that does not involve the exact root. Indeed, we implement the following adaptive multi-precision arithmetic scheme:

\vspace{-2mm}
\begin{equation}\label{eqdig5}
{\tt Digits } := \left[ \frac{\rho^2}{\rho-1} \times \left(- \log \left| \, \frac{f(x_{n})}{f(x_{n-1})} \right| \,+\, 2 \right) \,\right].
\end{equation}

Moreover, from (\ref{prop5}) we propose the following stopping criterion, instead of (\ref{stop1}):

\vspace{-2mm}
\begin{equation}\label{stop5}
     \left| \, \frac{f(x_{n})}{f(x_{n-1})} \right| <  10^{-\eta \: (\rho -1) / \rho} .
\end{equation}

%%%%%%%%%%%%%%%%%%%%%%%%%%%%%%%%%%%%%%%%%%%%%%%%%%%%%%%%%%%%%%%%%%%%%%%%%%%%%%%%%%%%%%%%%%%%%%%%%%%%%%%
%%%%%%%%%%%%%%%%%%%%%%%%%%%%%%%%%%%%%%%%%     Table 1     %%%%%%%%%%%%%%%%%%%%%%%%%%%%%%%%%%%%%%%%%%%%%

\vspace{-3mm}
\begin{table}[ht!!]
 \caption{\footnotesize Test functions, their roots and the initial points considered}\label{funs}
\vspace{-3mm}
{\footnotesize
  \begin{center}
   \begin{tabular}{lllc} \hline\hline  \\[-0.9em]
  \hspace{1.5cm} $f(x)$ &  \hspace{2cm}$\alpha$ & $x_0$ & $ \{x_{\scriptscriptstyle{-1}}\,,\, x_0 \}$ \\ [0.2em] \hline\\[-0.75em]
  $f_1(x)=x^3-3x^2+x-2$ & $2.893289196304497788906356$ & $2.5$ & $\{2.25, 2.60\}$ \\[0.1em]
  $f_2(x)=x^3+\cos x-2$ & $1.172577964753970012673333$ & $1.5$ & $\{1.50, 2.50\}$ \\[0.1em]
  $f_3(x)=2\sin x+1-x$  & $2.380061273139339017212548$ & $2.5$ & $\{1.00, 2.00\}$ \\[0.1em]
$f_4(x)=(x+1)\,e^{x-1}-1$ & $0.557145598997611416858672$& $1.0$ &$\{0.00, 0.75\}$ \\[0.1em]
$f_5(x)=e^{x^2+7x-30}-1$&                        $3.0$ & $2.94$& $\{2.90, 3.10\}$ \\[0.1em]
$f_6(x)=e^{-x}+\cos x.$ & $1.746139530408012417650703$ & $1.5$ & $\{1.60, 1.90\}$ \\[0.1em]
      $f_7(x)=x-3\ln x$ & $1.857183860207835336456981$ & $2.0$ & $\{1.00, 2.00\}$ \\ [0.3em]
    \hline \hline
  \end{tabular}
 \end{center}
 }
\end{table}

%
%
%%%%%%%%%%%%%%%%%%%%%%%%%%%%%%%%%%%%%%%%%%%%%%%%%%%%%%%%%%%%%%%%%%%%%%%%%%%%%%%%%%%%%%%%%%%%%%%%%%%
%%%%%%%%%%%%%%%%%%%%%%%%%%%%%%%%%%%%%%%%%%%%%%%%%%%%%%%%%%%%%%%%%%%%%%%%%%%%%%%%%%%%%%%%%%%%%%%%%%%
%
%

\vspace{-5mm}
\section{\large Iterative methods and numerical results}

\vspace{-2mm}
We consider in this section six iterative methods, $\phi_k, \, k=1\div 6$,  with local convergence order equal to $2, \,3, \,4, \,\frac{1+\sqrt{5}}{2},\, 1+\sqrt{2}$ and $1+\sqrt{3}$ respectively. The first three methods are one point iterative methods without memory known as Newton's method, Chebyshev's method \cite{GD5} and Schr\"{o}der's method \cite{sch}. The other three schemes are iterative methods with memory, namely the Secant method and two variants (see \cite{mgmn}). They are defined as

\vspace{-5mm}
\begin{eqnarray} \label{mi1}
\phi_1 (x_n) & = & x_n -\, u(x_n), \\[0.3em] \label{mi2}
\phi_2 (x_n) & = & \phi_1 (x_n) - \frac{1}{2} \,L(x_n)\, u(x_n), \\ \label{mi3}
\phi_3 (x_n) & = & \phi_2 (x_n) - \left( \frac{1}{2} \,L(x_n)^2\, - M(x_n) \right) \, u(x_n), \\ \label{mi4}
\phi_4 (x_n) & = & x_n -\, [x_{n-1},x_n]^{-1}_f \:f(x_n), \\[0.3em] \label{mi5}
\phi_5 (x_n) & = & \phi_4 (x_n) -\, [x_{n},\phi_4 (x_n)]^{-1}_f \:f(\phi_4 (x_n)), \\ [0.3em] \label{mi6}
\phi_6 (x_n) & = & \phi_4 (x_n) -\, [x_{n},2\phi_4 (x_n)-x_n]^{-1}_f \:f(\phi_4 (x_n)).
\end{eqnarray}

\vspace{-3mm}
\noindent where

\vspace{-3mm}
$$
 u(x) = \frac{f(x)}{f ' (x)}, \; \;\, L(x) = \frac{f^{\prime\prime}(x)}{f '(x)} \: u(x), \; \;\,  M(x) = \frac{f^{\prime\prime\prime}(x)}{3! \, f '(x)} \: u(x)^2 \; \; \:  {\rm and} \;  \; \; [x,y]^{-1}_f=\dfrac{y - x}{f(y) - f(x)} .
$$

%%
%%%%%
%%%%%%%%%%%%%%%%%%%%%%%%%%%%%%%%%%%%%%%%%%%%%%%%%%%%%%%%%%%%%%%%%%%%%%%%%%%%%%%%%%%%%%%%%%%%%%%
\vspace{1mm}
\noindent We point out that CLOC and PCLOC can be computed if $n \ge 1$, ACLOC if $n \ge 2$ and ECLOC if $n \ge 3$. If we have a method of higher order of convergence then multi-precision arithmetic is required and will be used where necessary to obtain many correct figures. In general, from guess points $ \, x_{-j}, \ldots, x_{\scriptscriptstyle -1},\,x_{\scriptscriptstyle 0}$ we obtain the admissible points $x_{\scriptscriptstyle 1},\ldots, x_{\scriptscriptstyle I}$. Notice that if we use an arithmetic with at most $\eta$ decimal digits,
with stopping criterion $|e_{\scriptscriptstyle I+1}|<10^{-\eta}$, then $x_{I+1}$ will not be considered since it would need a mantissa with higher length to hold all the correct decimals.

Hence, as $x_I$ is the best admissible point, we will take $n=I$ in the definitions of CLOC, ACLOC, ECLOC and PCLOC:

\vspace{-2mm}
\begin{defi} The computational values considered in all numerical experiments are

\vspace{-1mm}
\begin{equation}\label{TCLOC} {\overline \lambda} =
 \frac{\log \,\left| \, { e}_{\scriptscriptstyle I} \right|}{\log \, \left| \, { e}_{\scriptscriptstyle I-1} \right| } \, , \quad {\widehat \lambda} =
 \frac{\log \,\left| \, {\widehat e}_{\scriptscriptstyle I} \right|}{\log \, \left| \, {\widehat e}_{\scriptscriptstyle I-1} \right| } \, , \quad
{\widetilde \lambda} =
 \frac{\log \,\left| \, {\widetilde e}_{\scriptscriptstyle I} \right|}{\log \, \left| \, {\widetilde e}_{\scriptscriptstyle I-1} \right| } \quad \; {\rm and}\quad
\breve{\lambda} =\dfrac {\log |f({x}_{\scriptscriptstyle I})|}{\log |f(x_{\scriptscriptstyle{I-1}})|} \, .
\end{equation}

\end{defi}

\vspace{-1mm}
We have tested the preceding methods on seven functions using the Maple computer algebra system. We have computed the root of each function starting from the same initial approximation $\, x_0$ in (\ref{mi1})--(\ref{mi3}), and $ \{x_{-1}\,,\, x_0 \}$ in (\ref{mi4})--(\ref{mi6}). Depending on the computational order of convergence used, CLOC (\ref{CLOC}), ACLOC (\ref{ACLOC}), ECLOC (\ref{ECLOC}) or PCLOC (\ref{PCLOC}), we stop the iterative method when condition (\ref{stop1}), (\ref{stop3}), (\ref{stop2}) or (\ref{stop5}) is fulfilled. Note that in all cases  $\, \eta = 2200$.

\vspace{1mm}
\noindent The set of  test functions  presented here were previously considered in \cite{GS}. Table 1 shows the expression of each of these functions, their root with twenty five significant digits and their initial approximation. The latter is the same for all the iterative methods considered, considering one or two guess points depending on wether the algorithm works without or with memory.

%
%
%%%%%%%%%%%%%%%%%%%%%%%%%%%%%%%%%%%%%%%%%%%%%%%%%%%%%%%%%%%%%%%%%%%%%%%%%%%%%%%%%%%%%%%%%%%%%%%%%%%%%%%%
%%%%%%%%%%%%%%%%%%%%%%%%%%%%%%%%%%%%%%%%%%%%%%%%%%%%%%%%%%%%%%%%%%%%%%%%%%%%%%%%%%%%%%%%%%%%%%%%%%%%%%%%
%
%

\vspace{1mm}
For each method and function, we have applied the four techniques with adaptive multi-precision arithmetic (\ref{eqdig1}), (\ref{eqdig3}), (\ref{eqdig2}) and (\ref{eqdig5}). The number of necessary iterations to get the desired precision and the values of iterated points $x_1, \ldots, x_I$ are the same. Table 2 shows the number of iterations needed to compute the root. In addition, the last four columns show the interval with minimum and maximum error produced in the computation of the corresponding Computational Local Orders of Convergence (CLOC, ECLOC, ACLOC or PCLOC) for the seven test functions. For instance, considering the CLOC and Newton's method $\phi_1$, let us denote
$\mathcal{I}\, ( {\bar{\lambda}} )$ the resulting error interval obtained in the computation of the CLOC: $\,[ \min |\bar{\lambda}_k-\rho|, \:\max |\bar{\lambda}_k-\rho| ]$, for each function $f_k$, $k=1,\dots 7$.

%%%%%%%%%%%%%%%%%%%%%%%%%%%%%%%%%%%%%%%%%%%%%%%%%%%%%%%%%%%%%%%%%%%%%%%%%%%%%%%%%%%%%%%%%%%%%%%%%%%%%%%%
%%%%%%%%%%%%%%%%%%%%%%%%%%%%%%%%%%%%%%%%%%%%%%%      Table 2     %%%%%%%%%%%%%%%%%%%%%%%%%%%%%%%%%%%%%%%

\vspace{-3mm}
\begin{table}[ht!!]
 \caption{\small  Min-max interval for error bounds}\label{case7}
 \vspace{-6mm}
{\footnotesize
  \begin{center}
   \begin{tabular}{cccccccccccccc}
   \hline\hline \\ [-0.8em]
   &$f_1$&$f_2$&$f_3$&$f_4$&$f_5$&$f_6$&$f_7$&
    $\mathcal{I}\, ( \bar{\lambda} )$ & $\mathcal{I}\, ( \widetilde{\lambda} )$ & $\mathcal{I}\, ( \widehat{\lambda} )$ & $\mathcal{I}\, ( \breve{\lambda} )$\\[0.2em]\hline \\[-0.6em]
$\phi_1$&$12$&$11$&$10$&$11$&$12$&$10$&$11$& $[2.8e{-5},1.1e{-3}]$ & $[3.7e{-5},1.5e{-3}]$ & $[5.6e{-5},2.2e{-3}]$ & $[6.0e{-5},1.2e{-3}]$\\[0.2em]
$\phi_2$&$ 8$&$ 7$&$ 6$&$ 7$&$ 8$&$ 6$&$ 7$& $[8.9e{-5},3.3e{-3}]$ & $[1.0e{-4},6.0e{-3}]$ & $[1.4e{-4},9.9e{-3}]$ & $[2.1e{-4},4.5e{-3}]$\\[0.2em]
$\phi_3$&$ 6$&$ 6$&$ 5$&$ 6$&$ 6$&$ 5$&$ 5$& $[8.8e{-6},1.3e{-2}]$ & $[2.0e{-5},3.0e{-2}]$ & $[3.5e{-5},5.1e{-2}]$ & $[1.6e{-3},1.2e{-2}]$\\[0.2em]
$\phi_4$&$17$&$18$&$16$&$16$&$18$&$14$&$16$& $[8.1e{-6},5.8e{-4}]$ & $[7.9e{-6},6.8e{-4}]$ & $[1.2e{-5},9.4e{-4}]$ & $[3.2e{-5},5.5e{-3}]$\\[0.2em]
$\phi_5$&$ 9$&$ 9$&$ 9$&$ 8$&$10$&$ 7$&$ 8$& $[5.5e{-5},3.0e{-3}]$ & $[1.0e{-4},4.4e{-3}]$ & $[1.3e{-4},7.2e{-3}]$ & $[1.9e{-3},3.3e{-3}]$\\[0.2em]
$\phi_6$&$ 8$&$ 8$&$ 7$&$ 7$&$ 8$&$ 6$&$ 7$& $[3.6e{-5},3.7e{-3}]$ & $[1.4e{-4},1.6e{-2}]$ & $[1.6e{-4},1.7e{-2}]$ & $[2.8e{-4},4.5e{-3}]$\\[0.2em]
  \hline \hline
   \end{tabular}   			 	                 	
  \end{center}
}
 \end{table}

\vspace{-1mm}
From these numerical tests, we can conclude that the CLOC produces the best approximations of the theoretical order of convergence of an iterative method. However, the knowledge of the root is required. Conversely, as we can see in the definitions of ACLOC (\ref{ACLOC}), ECLOC (\ref{ECLOC}) and PCLOC (\ref{PCLOC}), these parameters have the advantage that they do not involve the expression of the root $\alpha$. Actually,  in real problems we want to approximate the root which it  is not know in advance. For practical purposes (see table \ref{case7} and Propositions  \ref{P41}, \ref{P51} and \ref{P61}) we recommend ECLOC since it presents the best approximation of the local order. Nevertheless, PCLOC is a good practical parameter in many cases because it requires less operations.

%
%
%%%%%%%%%%%%%%%%%%%%%%%%%%%%%%%%%%%%%%%%%%%%%%%%%%%%%%%%%%%%%%%%%%%%%%%%%%%%%%%%%%%%%%%%%%%%%%%%%%%
%%%%%%%%%%%%%%%%%%%%%%%%%%%%%%%%%%%%%%%%%%%%%%%%%%%%%%%%%%%%%%%%%%%%%%%%%%%%%%%%%%%%%%%%%%%%%%%%%%%
%
%

%%%%%%%%%%%%%%%%%%%%%%%%%%%%%%%%%%%%%%%%%%%%%%%%%%%%%%%%%%%%%%%%%%%%%%%%

\newpage

%%%%%%%%%%%%%%%%%%%%%%%%%%%%%%%%%%%%%%%%%%%%%%%%%%%%%%%%%%%%%%%%%%%%%%%%

\appendix
{\footnotesize
\section{\large Numerical computations}

In this section we show the results of the six iterative methods, $\phi_k, k=1\div 6,$ cited above applied to the test functions that appear in  table \ref{funs}. In each table we can find the necessary iteration number $I$ to get the required precision; the error $\Delta\bar{\lambda} = |\rho - \bar{\lambda}|$ made by the CLOC $\bar{\lambda}$; the error $\Delta\widehat{\lambda} = |\rho - \widehat{\lambda}|$ made by the ACLOC $\widehat{\lambda}$;
the error $\Delta\widetilde{\lambda}=|\rho - \widetilde{\lambda} |$ made by the ECLOC $\widetilde{\lambda}$ and the error $\Delta\breve{\lambda} = |\rho - \breve{\lambda}|$ made by the PCLOC $\breve{\lambda}$.

\vspace{1mm}
In tables \ref{case1}-\ref{case6} we give the details of the error intervals observed in the computations of CLOC, ACLOC, ECLOC and PCLOC shown in table \ref{case7}.

%%
%%%%%%%%%%%%%%
%%

\subsection{\normalsize Newton method -$\mathbf{\phi_1}$ method-}
The iterative Newton method defined by
\vspace{-2mm}
\begin{equation}
  \left\{
  \begin{array}{ll}
  {\rm given \ }x_0\,,\\ [0.4em]
   \phi_1 (x_{n}) =  \: x_{n} - \dfrac{f(x_{n})}{f^{\prime}(x_{n})},\quad n\geq0\,,
   \end{array}
  \right.
\end{equation}

\vspace{-2mm}
is used for the set of test functions (see table \ref{funs}). The errors $\Delta\bar{\lambda}$, $\Delta\widehat{\lambda}$, $\Delta\widetilde{\lambda}$ and $\Delta\breve{\lambda}$ obtained are shown in table \ref{case1}.

\vspace{-2mm}
%%%%%%%%%%%%%%%%%%%%%%%%%%%%%%%%%%%%%%%%  Table 3 %%%%%%%%%%%%%%%%%%%%%%%%%%%%%%%%%%%%%%%%%%%%%%%%%%%%%%%%%%
\begin{table}[ht!!]
{\footnotesize
 \caption{\footnotesize The errors $\Delta\bar{\lambda}$, $\Delta\widehat{\lambda}$, $\Delta\widetilde{\lambda}$ and $\Delta\breve{\lambda}$ for Newton's method}\label{case1} \vspace{-2mm}
  \begin{center}
   \begin{tabular}{ccccccccc}\hline\hline \\[-0.8em]
   & $I$ & $\Delta\bar{\lambda}$ & $\Delta\widehat{\lambda}$
         & $\Delta\widetilde{\lambda}$ & $\Delta\breve{\lambda}$ \\ \hline \\[-0.6em]
$f_1(x)$ & $12$ &$1.803\,e{-4}$&$3.607\,e{-4}$&$2.404\,e{-4}$&$1.086\,e{-3}$\\[0.1em]
$f_2(x)$ & $11$ &$2.790\,e{-5}$&$5.580\,e{-5}$&$3.720\,e{-5}$&$8.504\,e{-4}$\\[0.1em]
$f_3(x)$ & $10$ &$7.143\,e{-4}$&$1.430\,e{-3}$&$9.526\,e{-4}$&$1.220\,e{-3}$\\[0.1em]
$f_4(x)$ & $11$ &$2.723\,e{-4}$&$5.448\,e{-4}$&$3.632\,e{-4}$&$6.446\,e{-4}$\\[0.1em]
$f_5(x)$ & $12$ &$1.109\,e{-3}$&$2.215\,e{-3}$&$1.478\,e{-3}$&$4.018\,e{-4}$\\[0.1em]
$f_6(x)$ & $10$ &$1.040\,e{-3}$&$2.082\,e{-3}$&$1.387\,e{-3}$&$1.121\,e{-3}$\\[0.1em]
$f_7(x)$ & $11$ &$1.512\,e{-4}$&$3.025\,e{-4}$&$2.016\,e{-4}$&$6.032\,e{-5}$\\[0.1em]
    \hline \hline
  \end{tabular}
 \end{center}
 }
\end{table}

\vspace{-8mm}
%%%%%%%%%%%%%%
%%
\subsection{\normalsize Chebyshev method -$\mathbf{\phi_2}$ method-}
The iterative scheme of Chebyshev method is
\begin{equation}
 \left\{
  \begin{array}{llll}
                 {\rm given} & x_0\,,\\ [0.4em]
         y_{n} \,&=&\, \: x_{n} - u(x_{n})\,,\\[0.5em]
       \phi_2 (x_{n}) \,&=&\, \: y_{n} -\,\dfrac{1}{2}\, L(x_{n} )\; u(x_{n})\,,& n\geq0 \,.
  \end{array}
 \right.
\end{equation}
 where ${\displaystyle u(x) = {f(x)}/{f^{\prime}(x)}}$ and
       ${\displaystyle L(x) = \left({f^{\prime\prime}(x)}/{f^{\prime}(x)}\right) u(x)}.$
       The numerical results obtained for the set of test functions are listed in table \ref{case2}.

\vspace{-4mm}
%%%%%%%%%%%%%%%%%%%%%%%%%%%%%%%%%%%%%%%%  Table 4 %%%%%%%%%%%%%%%%%%%%%%%%%%%%%%%%%%%%%%%%%%%%%%%%%%%%%%%%%%
\begin{table}[ht!!]
{\footnotesize
\caption{\footnotesize The errors $\Delta\bar{\lambda}$, $\Delta\widehat{\lambda}$, $\Delta\widetilde{\lambda}$ and $\Delta\breve{\lambda}$ for Chebyshev's method}\label{case2} \vspace{-2mm}
  \begin{center}
   \begin{tabular}{cccccccc}\hline\hline \\[-0.8em]
   & $I$ & $\Delta\bar{\lambda}$ & $\Delta\widehat{\lambda}$
         & $\Delta\widetilde{\lambda}$ & $\Delta\breve{\lambda}$ \\ \hline \\[-0.6em]
$f_1(x)$ & $8$ &$2.077\,e{-4}$&$6.233\,e{-4}$&$3.739\,e{-4}$&$3.048\,e{-3}$\\[0.1em]
$f_2(x)$ & $7$ &$7.185\,e{-4}$&$2.154\,e{-3}$&$1.293\,e{-3}$&$2.148\,e{-3}$\\[0.1em]
$f_3(x)$ & $6$ &$1.949\,e{-3}$&$5.858\,e{-3}$&$3.511\,e{-3}$&$4.527\,e{-3}$\\[0.1em]
$f_4(x)$ & $7$ &$8.917\,e{-5}$&$1.437\,e{-4}$&$1.033\,e{-4}$&$2.109\,e{-4}$\\[0.1em]
$f_5(x)$ & $8$ &$3.318\,e{-3}$&$9.921\,e{-3}$&$5.965\,e{-3}$&$8.876\,e{-4}$\\[0.1em]
$f_6(x)$ & $6$ &$3.107\,e{-3}$&$9.350\,e{-3}$&$5.600\,e{-3}$&$3.612\,e{-3}$\\[0.1em]
$f_7(x)$ & $7$ &$2.017\,e{-4}$&$6.051\,e{-4}$&$3.630\,e{-4}$&$4.643\,e{-4}$\\[0.1em]
    \hline \hline
   \end{tabular}

  \end{center}
  }
 \end{table}

%%%%%%%%%%%%%%
%%

\subsection{\normalsize Schr\"{o}der method -$\mathbf{\phi_3}$ method-}

The iterative Schr\"{o}der method is defined by
\begin{equation}
 \left\{
  \begin{array}{llll}
                 {\rm given} & x_0\,,\\ [0.4em]
         y_{n}  \,&=&\, \: x_{n} - u(x_{n})\,,\\ [0.2em]
         z_{n}  \,&=&\, \: y_{n} - \dfrac{1}{2}\, L(x_{n} )\; u(x_{n})\,,\\ [0.6em]
         \phi_3 (x_{n})\,&=&\, \: z_{n} - \left(\dfrac{1}{2}\, \left(L(x_{n})\right)^2-M(x_{n})\right) \, u(x_{n}), &n\geq0 \,,
  \end{array}
  \right.
\end{equation}
where ${\displaystyle u(x) = \dfrac{f(x)}{f^{\prime}(x)}}$,
       ${\displaystyle L(x) = \dfrac{f^{\prime\prime}(x)} {f^{\prime}(x)}\,u(x)} \,$
and  ${\displaystyle M(x) = \dfrac{f^{\prime\prime\prime}(x)} {3!f^{\prime}(x)}\,\left(u(x)\right)^2\,.}$
The errors obtained in the approximations of $\rho$ by expressions (\ref{CLOC}), (\ref{ACLOC}), (\ref{ECLOC}) and (\ref{PCLOC}) for the set of test functions (see table \ref{funs}) are given in detail in table \ref{case3}.

%%%%%%%%%%%%%%%%%%%%%%%%%%%%%%%%%%%%%%%%  Table 5 %%%%%%%%%%%%%%%%%%%%%%%%%%%%%%%%%%%%%%%%%%%%%%%%%%%%%%%%%%
\begin{table}[ht!!]
{\footnotesize
\caption{\footnotesize The errors $\Delta\bar{\lambda}$, $\Delta\widehat{\lambda}$, $\Delta\widetilde{\lambda}$ and $\Delta\breve{\lambda}$ for Schr\"{o}der's method}\label{case3}     \vspace{-2mm}
  \begin{center}
   \begin{tabular}{cccccccc} \hline\hline\\[-0.8em]
   & $I$ & $\Delta\bar{\lambda}$ & $\Delta\widehat{\lambda}$
         & $\Delta\widetilde{\lambda}$ & $\Delta\breve{\lambda}$ \\ \hline \\[-0.6em]
$f_1(x)$ & $6$ &$8.809\,e{-6}$&$3.524\,e{-5}$  &$2.014\,e{-5}$&$1.218\,e{-2}$\\[0.1em]
$f_2(x)$ & $6$ &$1.261\,e{-3}$&$5.039\,e{-3}$  &$2.881\,e{-3}$&$2.091\,e{-3}$\\[0.1em]
$f_3(x)$ & $5$ &$2.142\,e{-3}$&$8.585\,e{-3}$  &$4.900\,e{-3}$&$6.033\,e{-3}$\\[0.1em]
$f_4(x)$ & $6$ &$2.635\,e{-4}$&$1.054\,e{-3}$  &$6.024\,e{-4}$&$1.567\,e{-3}$\\[0.1em]
$f_5(x)$ & $6$ &$1.299\,e{-2}$&$5.129\,e{-2}$  &$2.952\,e{-2}$&$2.814\,e{-3}$\\[0.1em]
$f_6(x)$ & $5$ &$3.778\,e{-3}$&$1.517\,e{-2}$  &$8.650\,e{-3}$&$4.521\,e{-3}$\\[0.1em]
$f_7(x)$ & $5$ &$6.252\,e{-5}$&$2.501\,e{-4}$  &$1.429\,e{-4}$&$3.027\,e{-3}$\\[0.1em]
    \hline \hline
   \end{tabular}
    \end{center}
  }
 \end{table}
%%
%%%%%%%%%%%%%%
%%

\subsection{\normalsize The Secant method -$\mathbf{\phi_4}$ method-}

The fourth iterative function considered $\phi_4$  is the Secant method defined by
\begin{equation}
 \left\{
  \begin{array}{llll}
                 {\rm given} & x_{-1}\,, & x_0\,,\\ [0.2em]
                 \phi_4 (x_{n}) &=& x_{n} -\, \dfrac{x_n-x_{n-1}}{f(x_n)-f(x_{n-1})}\,f(x_n)\,,\quad n\geq0 \,.
     \end{array}
  \right.
\end{equation}
The errors gotten in the approximated computation of $\rho$ by expressions (\ref{CLOC}), (\ref{ACLOC}), (\ref{ECLOC}) and (\ref{PCLOC}) for the set of test functions (see table \ref{funs}) in the secant method are displayed in table \ref{case4}.

%%%%%%%%%%%%%%%%%%%%%%%%%%%%%%%%%%%%%%%%  Table 6 %%%%%%%%%%%%%%%%%%%%%%%%%%%%%%%%%%%%%%%%%%%%%%%%%%%%%%%%%%
\begin{table}[ht!!]
{\footnotesize
  \caption{\footnotesize The errors $\Delta\bar{\lambda}$, $\Delta\widehat{\lambda}$, $\Delta\widetilde{\lambda}$ and $\Delta\breve{\lambda}$ for Secant method}\label{case4}
 \vspace{-2mm}
  \begin{center}
   \begin{tabular}{cccccccc}\hline\hline\\[-0.8em]
   & $I$ & $\Delta\bar{\lambda}$ & $\Delta\widehat{\lambda}$
         & $\Delta\widetilde{\lambda}$ & $\Delta\breve{\lambda}$ \\ \hline \\[-0.6em]
$f_1(x)$ & $17$ &$9.045\,e{-5}$&$1.466\,e{-4}$&$1.064\,e{-4}$&$5.448\,e{-4}$\\[0.1em]
$f_2(x)$ & $18$ &$8.112\,e{-6}$&$1.159\,e{-5}$&$7.925\,e{-6}$&$2.223\,e{-4}$\\[0.1em]
$f_3(x)$ & $16$ &$3.777\,e{-4}$&$6.100\,e{-4}$&$4.396\,e{-4}$&$6.448\,e{-4}$\\[0.1em]
$f_4(x)$ & $16$ &$1.090\,e{-4}$&$1.788\,e{-4}$&$1.321\,e{-4}$&$2.588\,e{-4}$\\[0.1em]
$f_5(x)$ & $18$ &$5.817\,e{-4}$&$9.408\,e{-4}$&$6.811\,e{-4}$&$2.107\,e{-4}$\\[0.1em]
$f_6(x)$ & $14$ &$5.110\,e{-4}$&$8.333\,e{-4}$&$6.098\,e{-4}$&$5.510\,e{-3}$\\[0.1em]
$f_7(x)$ & $15$ &$8.050\,e{-5}$&$1.295\,e{-4}$&$9.285\,e{-5}$&$3.187\,e{-5}$\\[0.1em]
    \hline \hline
   \end{tabular}
  \end{center}
  }
 \end{table}

\vspace{-2mm}
%%%%%%%%%%%%%%
%%

\subsection{\normalsize First variant of the Secant method -$\phi_5$ method-}

The iterative method $\phi_5$ is written by the following equations
\begin{equation}
 \left\{
  \begin{array}{llll}
                 {\rm given} & x_{-1}\,, & x_{0}\,,\\ [0.4em]
 y_{n} \,&=&\, \: \phi_4 (x_{n}),  \\ [0.4em]
 \phi_5 (x_{n}) \,&=&\, \: y_{n} -\, \dfrac{y_n-x_{n}}{f(y_n)-f(x_{n})}\,f(y_n)\,,& n\geq0 \,.
 \end{array}  \right.
\end{equation}

We apply this method to the set of test functions (table \ref{funs}) and the errors $\Delta\bar{\lambda}$, $\Delta\widehat{\lambda}$, $\Delta\widetilde{\lambda}$ and $\Delta\breve{\lambda}$ obtained are shown in table \ref{case5}.

\vspace{-4mm}
%%%%%%%%%%%%%%%%%%%%%%%%%%%%%%%%%%%%%%%%  Table 7 %%%%%%%%%%%%%%%%%%%%%%%%%%%%%%%%%%%%%%%%%%%%%%%%%%%%%%%%%%
\begin{table}[ht!!]
{\footnotesize
\caption{\footnotesize The errors $\Delta\bar{\lambda}$, $\Delta\widehat{\lambda}$, $\Delta\widetilde{\lambda}$ and $\Delta\breve{\lambda}$ for ${\phi_5}$ method}\label{case5}
    \vspace{-2mm}
  \begin{center}
   \begin{tabular}{cccccccc} \hline\hline\\[-0.8em]
   & $I$ & $\Delta\bar{\lambda}$ & $\Delta\widehat{\lambda}$ & $\Delta\widetilde{\lambda}$ & $\Delta\breve{\lambda}$ \\ \hline \\[-0.6em]
$f_1(x)$ & $9$ &$3.573\,e{-4}$&$8.632\,e{-4}$&$5.448\,e{-4}$&$2.152\,e{-3}$\\[0.1em]
$f_2(x)$ & $9$ &$5.517\,e{-5}$&$1.419\,e{-4}$&$1.026\,e{-4}$&$1.715\,e{-3}$\\[0.1em]
$f_3(x)$ & $9$ &$8.826\,e{-4}$&$2.135\,e{-3}$&$1.349\,e{-3}$&$1.507\,e{-3}$\\[0.1em]
$f_4(x)$ & $8$ &$5.266\,e{-4}$&$1.271\,e{-4}$&$7.997\,e{-4}$&$1.247\,e{-3}$\\[0.1em]
$f_5(x)$ &$10$ &$1.536\,e{-3}$&$3.702\,e{-3}$&$2.337\,e{-3}$&$5.563\,e{-4}$\\[0.1em]
$f_6(x)$ & $7$ &$3.033\,e{-3}$&$7.230\,e{-3}$&$4.375\,e{-3}$&$3.269\,e{-3}$\\[0.1em]
$f_7(x)$ & $8$ &$4.718\,e{-4}$&$1.136\,e{-3}$&$7.120\,e{-4}$&$1.876\,e{-4}$\\[0.1em]
    \hline \hline
  \end{tabular}
 \end{center}
 }
\end{table}

\vspace{-4mm}
%%%%%%%%%%%%%%
%%

\subsection{\normalsize Second variant of the Secant method -$\phi_6$ method-}

The iterative method $\phi_6$ is defined by the following iterative scheme
 \begin{equation}
 \left\{
  \begin{array}{llll}
                 {\rm given} & x_{-1}\,, & x_{0}\,,\\ [0.4em]
 y_{n} \,&=&\, \: \phi_4 (x_{n}),  \\ [0.4em]
\phi_6 (x_{n}) \,&=&\, \: y_{n} -\, \dfrac{2(y_{n}-x_{n})}{f(2y_{n}-x_{n})-f(x_{n})}\,f(y_{n})\,,& n\geq0 \,.
 \end{array}   \right.
\end{equation}

The errors $\Delta\bar{\lambda}$, $\Delta\widehat{\lambda}$, $\Delta\widetilde{\lambda}$ and $\Delta\breve{\lambda}$ of each sequence obtained by this method for each function of the set (see table \ref{funs}) are presented in table \ref{case6}.

\vspace{-4mm}
%%%%%%%%%%%%%%%%%%%%%%%%%%%%%% %%%%%%%%%  Table 8 %%%%%%%%%%%%%%%%%%%%%%%%%%%%%%%%%%%%%%%%%%%%%%%%%%%%%%%%%%
\begin{table}[ht!!]
{\footnotesize
  \caption{\footnotesize The errors $\Delta\bar{\lambda}$, $\Delta\widehat{\lambda}$, $\Delta\widetilde{\lambda}$ and $\Delta\breve{\lambda}$ for ${\phi_6}$ method}\label{case6}  \vspace{-2mm}
  \begin{center}
   \begin{tabular}{cccccccc} \hline\hline\\[-0.8em]
   & $I$ & $\Delta\bar{\lambda}$ & $\Delta\widehat{\lambda}$ & $\Delta\widetilde{\lambda}$ & $\Delta\breve{\lambda}$ \\ \hline \\[-0.6em]
$f_1(x)$ & $8$ &$4.641\,e{-4}$&$1.173\,e{-3}$&$6.377\,e{-4}$&$2.721\,e{-3}$\\[0.1em]
$f_2(x)$ & $8$ &$3.614\,e{-5}$&$1.570\,e{-4}$&$1.448\,e{-4}$&$1.385\,e{-3}$\\[0.1em]
$f_3(x)$ & $7$ &$2.626\,e{-3}$&$7.252\,e{-3}$&$4.481\,e{-3}$&$4.493\,e{-3}$\\[0.1em]
$f_4(x)$ & $7$ &$7.782\,e{-4}$&$1.477\,e{-3}$&$3.595\,e{-4}$&$1.705\,e{-3}$\\[0.1em]
$f_5(x)$ & $8$ &$3.405\,e{-3}$&$8.740\,e{-3}$&$4.913\,e{-3}$&$1.122\,e{-3}$\\[0.1em]
$f_6(x)$ & $6$ &$3.655\,e{-3}$&$1.672\,e{-2}$&$1.580\,e{-2}$&$4.021\,e{-3}$\\[0.1em]
$f_7(x)$ & $7$ &$6.321\,e{-4}$&$1.870\,e{-3}$&$1.262\,e{-3}$&$2.826\,e{-4}$\\[0.1em]
    \hline \hline
\end{tabular}
    \end{center}
  }
\end{table}

%
%
%%%%%%%%%%%%%%%%%%%%%%%%%%%%%%%%%%%%%%%%%%%%%%%%%%%%%%%%%%%%%%%%%%%%%%%%%%%%%%%%%%%%%%%%%%%%%%%%%%%%%%%
%%%%%%%%%%%%%%%%%%%%%%%%%%%%%%%%%%%%%%%%%%%%%%%%%%%%%%%%%%%%%%%%%%%%%%%%%%%%%%%%%%%%%%%%%%%%%%%%%%%%%%%
%
%
% ------------------------------------------------------------------------

}

% ------------------------------------------------------------------------

\vspace{-2mm}
\begin{thebibliography}{20}

{\footnotesize

\vspace{-1mm}
\bibitem{Wal} D.D. Wall,
The order of an iteration formula,
Math. Tables Aids Comput. 10 (1956) 167--168.

\vspace{-2.5mm}
\bibitem{Tor} L. Tornheim,
Convergence of multipoint iterative methods,
J. ACM 11 (1964) 210--220.

\vspace{-2.5mm}
\bibitem{OrRh} J.M. Ortega, W.C. Rheinboldt,
Iterative solution of nonlinear equations in several variables,
Academic Press, New York, 1970.

\vspace{-2.5mm}
\bibitem{WeFe} S. Weerakoon, T.G.I. Fernando,
A variant of Newton's method with accelerated third-order convergence,
Appl. Math. Lett. 13 (2000)  87--93.

\vspace{-2.5mm}
\bibitem{gn0} M. Grau, M. Noguera,
A variant of Cauchy's method with accelerated fifth-order convergence,
Appl. Math. Lett. 17 (2004) 509--517.

\vspace{-2.5mm}
\bibitem{chun} C. Chun,
Iterative methods improving Newton's method by the decomposition method,
Comput. Math. Appl. 50 (2005) 1559--1568.

\vspace{-2.5mm}
\bibitem{gng} M. Grau-S\'{a}nchez, M. Noguera, J.M. Guti\'{e}rrez,
On some computational orders of convergence,
Appl. Math. Lett. 23 (2010) 472--478.

\vspace{-2.5mm}
\bibitem{Aitk}  A. Aitken,
On Bernoulli's numerical solution of algebraic equations.
Proc. Roy. Soc. Edinburgh,  46 (1926) 289--305.

\vspace{-2.5mm}
\bibitem{Petk} M.S. Petkovi\'{c},
Remarks on "On a general class of multipoint root-finding methods of high computational efficiency",
SIAM J. Numer. Anal.  3 (2011) 1317--1319.

\vspace{-2.5mm}
\bibitem{DPP} J. D\u{z}uni\'{c}, M.S. Petkovi\'{c}, L.D. Petkovi\'{c},
Three-point methods with and without memory for solving nonlinear equations,
Appl. Math. Comput. (2011), doi:10.1016/j.amc.2011.10.057

\vspace{-2.5mm}
\bibitem{GD5} M. Grau, J. L. D\'{\i}az-Barrero,
An improvement of the Euler-Chebyshev iterative method,
J. Math. Anal. Appl.  315 (2006) 1--7.

\vspace{-2.5mm}
\bibitem{sch}  E. Schr\"{o}der,
\"{U}ber unendlich viele Algorithmen zur Aufl\"{o}sung der Gleichungen,
Math. Ann. 2 (1870) 317--365.

\vspace{-2.5mm}
\bibitem{mgmn} M. Grau-S\'{a}nchez, M. Noguera,
A technique to choose the most efficient method between secant method and some variants,
Appl. Math. Comput. 218 (2012) 6415--6426.


\vspace{-2.5mm}
\bibitem{GS} M. Grau-S\'{a}nchez,
Improvement of the efficiency of some three-step iterative like-Newton methods,
Numer. Math. 107 (2007) 131--146.

}
%
%
\end{thebibliography}
\end{document}